\documentclass[10pt]{amsart}
\usepackage{mathptmx}
\usepackage{amsmath}
\usepackage{amssymb}
\usepackage{array}
\usepackage{geometry}
\usepackage[bookmarks=true,colorlinks=true, pdfstartview=FitV, linkcolor=black, citecolor=blue, urlcolor=black]{hyperref}

\usepackage{color}
\definecolor{DarkRed}{rgb}{0.55,.00,0.2}
\definecolor{DarkGrey}{rgb}{0.35,.35,0.35}

\theoremstyle{definition}

\theoremstyle{remark}

\numberwithin{equation}{section}

\begin{document}

\title{ Lebedev's type   index transforms with the squares of the  associated Legendre functions}

\author{S. Yakubovich}
\address{Department of Mathematics, Faculty of Sciences,  University of Porto,  Campo Alegre str.,  687; 4169-007 Porto,  Portugal}
\email{ syakubov@fc.up.pt}

\keywords{Index Transforms,  Associated Legendre functions,   modified Bessel functions, Fourier transform,   Mellin transform, Mehler-Fock transform, Boundary  value problem}
\subjclass[2000]{  44A15, 33C10, 33C45, 44A05
}

\date{\today}
\maketitle

\markboth{\rm \centerline{ S.  Yakubovich}}{}
\markright{\rm \centerline{Index transforms with the squares and products of the associated Legendre functions}}

\begin{abstract}   The classical Lebedev index transform (1967), involving squares and products   of the Legendre functions is generalized on the associated Legendre functions of an arbitrary order.     Mapping properties  are investigated  in the Lebesgue  spaces.  Inversion formulas  are proved.    As an interesting application,  a solution to the boundary value problem for a  third  order partial differential equation is obtained.  
\end{abstract}

\section{Introduction and preliminary results}

In 1967 N.N. Lebedev  proved (cf. \cite{square}) that the  antiderivative of an  arbitrary function $f$ defined on the interval $(1,\infty)$ and satisfying integrability conditions $f \in  L_1((1,a); (x-1)^{-1/2} dx),\  f \in L_1 ((a,\infty); \log x\  dx)$ for some $ a >1$ can be expanded in terms of the following repeated integral for $x \in (1,\infty)$

$$\int_1^x f(y)dy = 2(x^2-1)^{1/2} \int_0^\infty \tau \tanh(\pi\tau) \left[ P_{-1/2+i\tau} (x)\right]^2  $$

$$\times \int_1^\infty (y^2-1)^{1/2}   P_{-1/2+i\tau} (y) \left[  Q_{-1/2+i\tau} (y) +  Q_{-1/2- i\tau} (y)\right] f(y) dy d\tau,\eqno(1.1)$$ 
where $P_\nu(z), Q_\nu(z)$ are Legendre functions of the first and second kind, respectively,  (see \cite{erd}, Vol. I, \cite{vir}). This expansion generates a pair of the direct and inverse  index transforms \cite{yak}, namely,

$$F(\tau)=  \int_1^\infty (y^2-1)^{1/2}   P_{-1/2+i\tau} (y) \left[  Q_{-1/2+i\tau} (y) +  Q_{-1/2- i\tau} (y)\right] f(y) dy\eqno(1.2)$$
for all $\tau >0$,

$$f(x)=  2 {d\over dx}  \int_0^\infty \tau \tanh(\pi\tau)  \left[ P_{-1/2+i\tau} (x)\right]^2  (x^2-1)^{1/2} F(\tau) d\tau\eqno(1.3)$$
for almost all $x>0$.   The main goal of the present paper is to extend Lebedev's  transforms (1.2), (1.3),   considering Lebedev-like kernels, which  involve products of the associated Legendre functions of the first and second kind of an arbitrary complex order $\mu$, correspondingly,   $P^\mu_\nu(z), Q^\mu_\nu(z)$  \cite{vir}. We will study them mapping properties, prove inversion theorems and apply to solve the boundary value  problem for a  higher order PDE.  Precisely, we will consider the following operator of the index transform

$$(F_\mu f) (\tau)=  \sqrt\pi\  \Gamma\left({1\over 2} + i\tau -\mu\right) \Gamma\left({1\over 2} - i\tau -\mu\right)
 \int_0^\infty \left[ P^\mu_{-1/2+i\tau} \left(\sqrt{1+y\over y}\right)\right]^2 f(y) dy,\  \tau \in \mathbb{R},\eqno(1.4)$$
 and its adjoint one

 $$(G_\mu g) (x)=  \sqrt\pi \int_{-\infty}^\infty \Gamma\left({1\over 2} + i\tau -\mu\right) \Gamma\left({1\over 2} - i\tau -\mu\right)  \left[ P^\mu_{-1/2+i\tau} \left(\sqrt{1+x\over x}\right)\right]^2 g(\tau) d\tau, \  x \in \mathbb{R}_+,\eqno(1.5)$$
 where $\Gamma(z)$ is Euler's gamma-function $\mu \in \mathbb{C}$, $i$ is the imaginary unit  and the integration in (1.5) is realized with respect to the lower index of the associated Legendre function of the first kind.   We note that the kernel $P^\mu_{-1/2+i\tau} (x), \ x >1$ corresponds to the well-known generalized Mehler-Fock transform \cite{yak}. Denoting the kernel of (1.4), (1.5) by

 $$\Phi_\tau(x)=  \sqrt\pi \  \Gamma\left({1\over 2} + i\tau -\mu\right) \Gamma\left({1\over 2} - i\tau -\mu\right) \left[ P^\mu_{-1/2+i\tau} \left(\sqrt{1+x\over x}\right)\right]^2,\eqno(1.6)$$
we will find for further use its representation in terms of Fourier cosine transform \cite{tit} and deduce an ordinary differential equation with polynomial coefficients, whose solution is $\Phi_\tau(x)$, employing  the so-called method of the Mellin-Barnes integrals, which is already being successfully applied by the author for other index transforms.  In fact, appealing to relation (8.4.41.47) in \cite{prud}, Vol. III, we find the following Mellin-Barnes integral representation for the kernel (1.6), namely

$$ \Phi_\tau(x) = {1\over 2\pi i} \int_{\gamma-i\infty}^{\gamma +i\infty} \frac {\Gamma(s+1/2+  i\tau)\Gamma(s+1/2 -i\tau) \Gamma(1/2+s)\Gamma(-\mu-s) }{\Gamma(1+s) \Gamma (1+s-\mu) } x^{-s} ds, \ x >0,\eqno(1.7)$$
where $\gamma$ is taken from the interval $(-1/2, - {\rm Re} \mu)$. The absolute convergence of the integral (1.7)   follows immediately from  the Stirling asymptotic formula for the gamma- function \cite{erd}, Vol. I, because for all $\tau \in \mathbb{R}$
$$\frac {\Gamma(s+1/2+  i\tau)\Gamma(s+1/2 -i\tau) \Gamma(1/2+s)\Gamma(-\mu-s) }{\Gamma(1+s) \Gamma (1+s-\mu) }   =  O \left(e^{-\pi |s|} |s|^ {-{\rm Re} \mu -3/2}\right),\  |s| \to \infty.\eqno(1.8)$$
Moreover, it can be differentiated under the integral sign any number of times due to the absolute and uniform convergence by $x \ge x_0 >0$.   We have

{\bf Lemma 1}. {\it Let $x \ge 1,  \tau \in \mathbb{R},\  {\rm Re}\mu < 1/2$. Then the kernel $(1.6)$ has the following representation in terms of  Fourier cosine transform of the second kind Legendre function, namely} 
$$ \Phi_\tau(x)  =  2 \sqrt{x\over \pi}  \int_0^\infty  \cos(\tau u)  Q_{-1/2-\mu} \left( 2 x \cosh^2(u/2)+1\right) du.\eqno(1.9) $$

\begin{proof}     In fact, appealing  to the reciprocal formulae via the Fourier cosine transform (cf. formula (1.104) in \cite{yak}) 
$$\int_0^\infty  \Gamma\left(s +1/2 + i\tau\right)  \Gamma\left(s +1/2 - i\tau \right)  \cos( \tau y) d\tau
= {\pi\over 2^{2s+1}}  {\Gamma(2s+1) \over \cosh^{2s+1}(y/2)},\ {\rm Re}\ s > - {1\over 2},\eqno(1.10)$$
$$  \Gamma\left(s +1/2+ i\tau \right)  \Gamma\left(s +1/2 - i\tau\right)  
=   { \Gamma(2s+1)  \over 2^{2s}}  \int_0^\infty   {\cos(\tau y)  \over \cosh^{2s+1} (y/2)} \ dy,\eqno(1.11)$$ 
we replace  the gamma-product $ \Gamma\left(s +1/2+ i\tau \right)  \Gamma\left(s+1/2 - i\tau\right)$ in the integral (1.9)  by its integral representation (1.11) and change the order of integration via Fubini's theorem.  Then, employing the duplication formula for the gamma- function \cite{erd}, Vol. I,   we derive 
$$ \Phi_\tau(x)  = {1\over 2 \pi\sqrt \pi  i} \int_0^\infty   {\cos(\tau u) \over \cosh(u/2)}  \int_{\gamma-i\infty}^{\gamma +i\infty} \frac {\left[\Gamma(s+ 1/2)\right]^2 \Gamma(-\mu-s)}{ \Gamma (1+s-\mu) } \left( x \cosh^2(u/2) \right)^{-s}  ds du  $$
$$=  2 \sqrt{x\over \pi}  \int_0^\infty  \cos(\tau u)  Q_{-1/2-\mu} \left( 2 x \cosh^2(u/2)+1\right) du,$$
where the  inner integral with respect to $s$ is calculated via Slater's theorem as a sum of residues at simple right-hand side poles of the gamma-function $\Gamma(-\mu-s)$ and relation (7.3.1.71) in \cite{prud}, Vol. III.   Lemma 1 is proved. 

\end{proof}

{\bf Lemma 2}.    {\it For each $\tau \in \mathbb{R}$ the function  $ \Phi_\tau(x)$ given by formula $(1.6)$  is  a fundamental solution of the following third order differential equation with polynomial   coefficients}
$$2x^3(1+x)  {d^3 \Phi_{\tau} \over dx^3} + 3x(1+2x) {d^2 \Phi_\tau \over dx^2}  +  x \left( 2x(1-\mu^2) + 2\tau^2+ {1\over 2} \right)  {d\Phi_\tau \over dx}  -  \left( {1\over 4} + \tau^2\right) \Phi_\tau = 0,\ x >0. \eqno(1.12)$$

\begin{proof}  As it was mentioned above, the asymptotic behavior (1.8) of the integrand in (1.7) permits a differentiation under the integral sign any number of times.  Hence employing the reduction formula for the gamma- function \cite{erd}, Vol. I , we derive 

$$\left( x {d\over dx} \right)^2  \Phi_{\tau} =   {1\over 2\pi i} \int_{\gamma-i\infty}^{\gamma +i\infty} \frac { s^2 \Gamma(s+1/2+  i\tau)\Gamma(s+1/2 -i\tau) \Gamma(1/2+s)\Gamma(-\mu-s) }{\Gamma(1+s) \Gamma (1+s-\mu) } x^{-s} ds$$

$$= {1\over 2\pi i} \int_{\gamma-i\infty}^{\gamma +i\infty} \frac {  \Gamma(s+3/2+  i\tau)\Gamma(s+3/2 -i\tau) \Gamma(1/2+s)\Gamma(-\mu-s) }{\Gamma(1+s) \Gamma (1+s-\mu) } x^{-s} ds -  \left({1\over 4}+ \tau^2\right)  \Phi_{\tau}  + x {d  \Phi_{\tau} \over dx}.\eqno(1.13)$$
Meanwhile, with a simple change of variable

$${1\over 2\pi i} \int_{\gamma-i\infty}^{\gamma +i\infty} \frac {  \Gamma(s+3/2+  i\tau)\Gamma(s+3/2 -i\tau) \Gamma(1/2+s)\Gamma(-\mu-s) }{\Gamma(1+s) \Gamma (1+s-\mu) } x^{-s} ds$$

$$= {1\over 2\pi i} \int_{1+\gamma -i\infty}^{1+\gamma +i\infty} \frac {  \Gamma(s+1/2+  i\tau)\Gamma(s+1/2 -i\tau) \Gamma(s-1/2)\Gamma(1-\mu-s) }{\Gamma(s) \Gamma (s-\mu) } x^{1-s} ds$$

$$= {1\over 2\pi i} \int_{1+\gamma -i\infty}^{1+\gamma +i\infty} \frac { s(s-\mu) (-\mu-s)  \Gamma(s+1/2+  i\tau)\Gamma(s+1/2 -i\tau) \Gamma(s+1/2)\Gamma(-\mu-s) }{(s-1/2)\Gamma(1+s) \Gamma (1+s-\mu) } x^{1-s} ds$$

$$= -  {1\over 2\pi i} \int_{1+\gamma -i\infty}^{1+\gamma +i\infty} \frac {(s^2-\mu^2)  \Gamma(s+1/2+  i\tau)\Gamma(s+1/2 -i\tau) \Gamma(s+1/2)\Gamma(-\mu-s) }{\Gamma(1+s) \Gamma (1+s-\mu) } x^{1-s} ds$$

$$-   {1\over 4\pi i} \int_{1+\gamma -i\infty}^{1+\gamma +i\infty} \frac {(s^2-\mu^2)  \Gamma(s+1/2+  i\tau)\Gamma(s+1/2 -i\tau) \Gamma(s+1/2)\Gamma(-\mu-s) }{(s-1/2) \Gamma(1+s) \Gamma (1+s-\mu) } x^{1-s} ds.$$
Therefore, moving the contour of  two latter integrals to the left via Cauchy's theorem, dividing by $\sqrt x$, differentiating again and using (1.13), we obtain

$${d\over dx}  \left[ {2\over \sqrt x} \left[ \left( x {d\over dx} \right)^2  \Phi_{\tau} +  \left({1\over 4}+ \tau^2\right)  \Phi_{\tau}  -  x {d  \Phi_{\tau} \over dx} + x \left( x {d\over dx} \right)^2  \Phi_{\tau} - x\mu^2  \Phi_{\tau} \right] \right]$$

$$= {1\over \sqrt x} \left[  \left( x {d\over dx} \right)^2  \Phi_{\tau}  -\mu^2  \Phi_{\tau} \right].$$
This is equivalent to the following operator equation

$$2(1+x)   \left( x {d\over dx} \right)^3  \Phi_{\tau} - 3 \left( x {d\over dx} \right)^2  \Phi_{\tau} + \left({3\over 2} +2(\tau^2-x\mu^2)\right) \left( x {d\over dx} \right)  \Phi_{\tau}  -  \left( {1\over 4} + \tau^2\right) \Phi_\tau = 0,\ x >0,$$
which drives to (1.12), fulfilling simple differentiation.   Lemma 2 is proved. 

\end{proof}

\section {Boundedness  and inversion properties for the index transform (1.4)}

In this section we will investigate the mapping properties of the index transform (1.4), involving  the Mellin transform technique developed in \cite{yal}.   Precisely, the Mellin transform is defined, for instance, in  $L_{\nu, p}(\mathbb{R}_+),\ 1 \le  p \le 2$ (see details in \cite{tit}) by the integral  
$$f^*(s)= \int_0^\infty f(x) x^{s-1} dx,\eqno(2.1)$$
 being convergent  in mean with respect to the norm in $L_q(\nu- i\infty, \nu + i\infty),\ \nu \in \mathbb{R}, \   q=p/(p-1)$.   Moreover, the  Parseval equality holds for $f \in L_{\nu, p}(\mathbb{R}_+),\  g \in L_{1-\nu, q}(\mathbb{R}_+)$
$$\int_0^\infty f(x) g(x) dx= {1\over 2\pi i} \int_{\nu- i\infty}^{\nu+i\infty} f^*(s) g^*(1-s) ds.\eqno(2.2)$$
The inverse Mellin transform is given accordingly
 $$f(x)= {1\over 2\pi i}  \int_{\nu- i\infty}^{\nu+i\infty} f^*(s)  x^{-s} ds,\eqno(2.3)$$
where the integral converges in mean with respect to the norm  in   $L_{\nu, p}(\mathbb{R}_+)$
$$||f||_{\nu,p} = \left( \int_0^\infty  |f(x)|^p x^{\nu p-1} dx\right)^{1/p}.\eqno(2.4)$$
In particular, letting $\nu= 1/p$ we get the usual space $L_1(\mathbb{R}_+)$.  Further, denoting by $C_b(\mathbb{R}_+)$ the space of bounded continuous functions, we establish

{\bf Theorem 1.}   {\it Let ${\rm Re}\mu < 1/2$. The index transform  $(1.4)$  is well-defined as a  bounded operator $F_\mu: L_{1-\nu, 1} \left(\mathbb{R}_+\right) \to C_b (\mathbb{R}),\  \nu \in (-1/2,\ - {\rm Re}\mu)$ and the following norm inequality takes place
$$||F_\mu f||_{C_b(\mathbb{R})} \equiv \sup_{\tau \in \mathbb{R}} | (F_\mu f)(\tau)| \le C_{\mu,\nu}  ||f||_{1-\nu,1},\eqno(2.5)$$
where

$$C_{\mu,\nu} =  {2^{2\nu -1} \over  \pi\sqrt \pi } B\left(\nu + {1\over 2},\  \nu + {1\over 2} \right) \int_{\nu-i\infty}^{\nu +i\infty} \left| \frac {\left[\Gamma(s+ 1/2)\right]^2 \Gamma(-\mu-s)}{ \Gamma (1+s-\mu) }  ds\right| \eqno(2.6)$$
and $B(a,b)$ is the Euler beta-function \cite{erd}, Vol. 1.  Moreover,  $(F_\mu f)(\tau) \to 0,\  |\tau| \to \infty.$  Besides, if, in addition,  $ f \in L_{1-\nu, p},\ 1 < p\le 2$ and $ \nu \in \left(- 1/(2q), \   \hbox{min} (0, - {\rm Re} \mu)\right),\ q= p/(p-1)$, then  

$$(F_\mu f)(\tau) =   {\sqrt{\pi} \over\cosh(\pi\tau)} \int_0^\infty   K_{i\tau} \left({1\over \sqrt x}\right) \left[  I_{i\tau} \left({1\over \sqrt x}\right) +  I_{- i\tau} \left({1\over \sqrt x}\right)\right] \varphi(x) {dx\over x},\eqno(2.7)$$
where 

$$\varphi(x)=  {1\over 2\pi i} \int_{\nu+1/2 -i\infty}^{\nu+1/2  +i\infty} \frac { \Gamma(1/2 -\mu-s)\Gamma(s) \Gamma(1-s)}{\Gamma(s+1/2) \Gamma(1/2 -s) \Gamma (1/2+s-\mu) } f^*(3/2 -s) x^{-s} ds\eqno(2.8)$$
and  integrals $(2.7),\ (2.8)$ converge absolutely.}

\begin{proof}  Indeed, following ideas elaborated in the proof of Lemma 1 and calculating an elementary integral with the hyperbolic function (see relation (2.4.4.4) in \cite{prud}, Vol. I), we have 

$$ \left|(F_\mu f)(\tau) \right| \le \int_0^\infty \left| \Phi_\tau(y)\right|  |f(y)| dy   \le {1\over 2 \pi\sqrt \pi } \int_0^\infty   {du  \over \left[\cosh(u/2)\right]^{2\gamma +1} } $$

$$\times  \int_{\nu-i\infty}^{\nu +i\infty} \left| \frac {\left[\Gamma(s+ 1/2)\right]^2 \Gamma(-\mu-s)}{ \Gamma (1+s-\mu) }  ds\right|  \int_0^\infty  |f(y)| y^{-\nu } dy = C_{\mu,\nu} ||f||_{1-\nu,1}. \eqno(2.9)$$
This proves (2.5).   Furthermore, Fubini's theorem and definition of the Mellin transform (2.1) allow us to write the composition representation for index transform (1.4) in terms of the Fourier cosine transform, namely,

$$(F_\mu f)(\tau) = {1\over 2 \pi\sqrt \pi  i} \int_0^\infty   {\cos(\tau u) \over \cosh(u/2)}  \int_{\nu-i\infty}^{\nu +i\infty} \frac {\left[\Gamma(s+ 1/2)\right]^2 \Gamma(-\mu-s)}{ \Gamma (1+s-\mu) } f^*(1-s) \left( \cosh^2(u/2) \right)^{-s}  ds du, $$
where

$$\Psi(u)=  \int_{\nu-i\infty}^{\nu +i\infty} \frac {\left[\Gamma(s+ 1/2)\right]^2 \Gamma(-\mu-s)}{ \Gamma (1+s-\mu) } f^*(1-s) \left( \cosh^2(u/2) \right)^{-s-1}  ds  \in L_1\left(\mathbb{R }_+;\ du\right).$$
Hence it tends to zero, when $|\tau| \to \infty$ via the Riemann-Lebesgue lemma.  On the other hand, returning to the Mellin-Barnes representation (1.7), appealing to the Parseval equality (2.2) and making a simple substitution, we deduce the formula 

$$ (F_\mu f)(\tau) =  {1\over 2\pi i} \int_{\nu+1/2-i\infty}^{\nu+1/2 +i\infty} \frac {\Gamma(s+  i\tau)\Gamma(s -i\tau) \Gamma(s)\Gamma(1/2 -\mu-s) }{\Gamma(1/2+s) \Gamma (1/2+s-\mu) } f^*(3/2-s) ds.\eqno(2.10)$$

Meanwhile,   taking the Mellin-Barnes representation for the product of modified Bessel functions of the third kind (cf. relation (8.4.23.24) in \cite{prud}, Vol. III)

$${\sqrt{\pi} \over x\cosh(\pi\tau)} K_{i\tau} \left({1\over \sqrt x}\right) \left[  I_{i\tau} \left({1\over \sqrt x}\right) +  I_{- i\tau} \left({1\over \sqrt x} \right)\right] $$

$$={1\over 2\pi i} \int_{\gamma -i\infty}^{\gamma  +i\infty} \frac { \Gamma( s-1/2) }{\Gamma(s)} \   \Gamma(1+  i\tau-s)\Gamma(1 -i\tau-s) x^{-s} ds,\ x >0,$$
where $1/2 < \gamma < 1$ and using again the Parseval equality (2.2), we find that (2.10) becomes the Lebedev index transform with the product of the modified Bessel functions \cite{square} given by formula (2.7),  where $\varphi(x)$ is  defined by integral (2.8). Its  convergence for each $x >0$  is absolute due to the estimate with the aid of H\"older's  inequality 

$$ \int_{\nu+1/2 -i\infty}^{\nu+1/2  +i\infty} \left|\frac { \Gamma(1/2 -\mu-s)\Gamma(s ) \Gamma(1 -s)}{\Gamma(s+1/2) \Gamma(1/2-s) \Gamma (1/2+s-\mu) } f^*(3/2-s) x^{-s} ds\right| \le x^{-\nu-1/2} \left( \int_{\nu -i\infty}^{\nu  +i\infty} \left|f^*(1-s)\right|^q |ds|\right)^{1/q}$$

$$\times \left( \int_{\nu -i\infty}^{\nu  +i\infty} \left|\frac { \Gamma( -\mu-s)\Gamma(s+1/2) \Gamma(1/2 -s)}{\Gamma(s+1) \Gamma(-s) \Gamma (1+s-\mu) } \right|^p |ds|\right)^{1/p}  < \infty,\ q= {p\over p-1}.$$
The convergence of the latter integral by $s$ is justified,  recalling  the Stirling formula for the asymptotic behavior of the gamma-function, which gives

$$ \frac { \Gamma(-\mu-s)\Gamma(s+1/2) \Gamma(1/2 -s)}{\Gamma(s+1) \Gamma(-s) \Gamma (1+s-\mu) } 
= O\left( |s|^{-2\nu-1}\right),\ |s| \to \infty$$
and $\nu$ is chosen from the interval  $\left(- 1/(2q), \   \hbox{min} (0, - {\rm Re} \mu)\right)$.  Moreover, we find that $\varphi (x)= O(x^{-\nu-1/2}),\ x >0.$ This asymptotic behavior is used to establish the absolute convergence of the integral (2.7). In fact, employing the asymptotic formulae for the modified Bessel functions \cite{erd}, Vol. 2 for fixed $\tau \in \mathbb{R}$

$$K_{i\tau} \left(\sqrt {x}\right) \left[  I_{i\tau} \left(\sqrt { x}\right) +  I_{- i\tau} \left(\sqrt { x}\right)\right] = O(\log x),\  x\to 0+, $$

$$ K_{i\tau} \left(\sqrt {x}\right) \left[  I_{i\tau} \left(\sqrt { x}\right) +  I_{- i\tau} \left(\sqrt { x}\right)\right] = O\left( x^{-1/2}\right),\ x \to \infty,$$
we derive via elementary substitution  

$$\int_0^\infty   \left| K_{i\tau} \left({1\over \sqrt x}\right) \left[  I_{i\tau} \left({1\over \sqrt x}\right) +  I_{- i\tau} \left({1\over \sqrt x}\right)\right] \varphi(x)\right| {dx\over x} = \int_0^1 O\left ( x^{\nu-1/2} \log x \right) dx $$

$$ +   \int_1^\infty  O\left ( x^{\nu-1}  \right) dx < \infty,\  \nu \in \left(- 1/(2q), \   \hbox{min} (0, - {\rm Re} \mu)\right).$$ 
Theorem 1 is proved. 

\end{proof}

Writing  (2.7) in the form

$$(F_\mu f)(\tau) =   {2 \sqrt{\pi} \over\cosh(\pi\tau)} \int_0^\infty   K_{i\tau} \left(x\right) \left[  I_{i\tau} \left(x\right) +  I_{- i\tau} \left( x\right)\right] \varphi(1/x^2){ dx\over x} ,\eqno(2.11)$$
we will appeal to the Lebedev expansion theorem in \cite{square},  which implies the following representation of the antiderivative 

$$ \int_{x}^{\infty}  \varphi\left( {1\over y^2}\right) {dy\over y}  = {1\over \pi^2\sqrt \pi} \int_0^\infty \tau \sinh(2\pi\tau) K_{i\tau}^2(x) (F_\mu f)(\tau) d\tau,\ x >0,\eqno(2.12)$$
which holds under condition $\varphi(1/x^2) x^{-1} \in L_1\left((0,1); \ x^{-1/2} dx\right) \cap  L_1\left((1,\infty); \  x^{1/2} dx\right)$. By straightforward substitutions we see that this condition is equivalent to  (cf. (2.4)) $\varphi \in L_{1/4,1}(1,\infty) \cap L_{-1/4,1}(0,1)$.   On the other hand, letting in (2.12) $1/\sqrt x$ instead of $x$ and changing the variable in its left-hand side,  we obtain 

$$ \int_0^{ x} \varphi\left( y\right) {dy\over y} = {2\over \pi^2\sqrt \pi} \int_0^\infty \tau \sinh(2\pi\tau) K_{i\tau}^2\left({1\over \sqrt x}\right)  (F_\mu f)(\tau)  d\tau.$$
Then, appealing to relation (8.4.23.28) in \cite{prud}, Vol. III and  differentiating two sides of the latter equality with respect to $x$, we find

$$\varphi(x) = {x \over 2\pi^3 i} {d\over dx}  \int_0^\infty \tau \sinh(2\pi\tau)  (F_\mu f)(\tau) 
\int_{\gamma -i\infty}^{\gamma   +i\infty} \frac { \Gamma(-s)\Gamma(i\tau-s) \Gamma(-i\tau -s)}{\Gamma(1/2 -s)} x^{-s} ds d\tau,\  \gamma <0.\eqno(2.13)$$
Our goal now is to motivate the differentiation under integral sign in the right-hand side of (2.13).   To do this we use (1.11) to get the uniform estimate of the product of gamma-functions 

$$ \left| \Gamma\left(i\tau-s \right)  \Gamma\left(- i\tau-s\right)\right| \le   
   2^{{\rm Re} s +1} \left|\Gamma(-2s)\right|  \int_0^\infty   \cosh^{2{\rm Re} s}  \left({y\over 2}\right) \ dy, \  {\rm Re} s < 0.$$ 
Consequently, assuming the integrability condition $F_\mu \in L_1(\mathbb{R};  |\tau| e^{2\pi|\tau|} d\tau)$, the repeated integral in (2.13) can be estimated as follows 

$$ \int_0^\infty \tau \sinh(2\pi\tau)  \left|(F_\mu f)(\tau)\right| 
\int_{\gamma -i\infty}^{\gamma  +i\infty} \left| \frac { \Gamma(-s)\Gamma(i\tau-s) \Gamma(-i\tau -s)}{\Gamma(1/2 -s)} x^{-s} ds\right| d\tau$$

$$ \le C x^{- \gamma}  ||F_\mu||_{ L_1(\mathbb{R}_+; |\tau | e^{2\pi \tau} d\tau)} \int_{\gamma -i\infty}^{\gamma  +i\infty} \left| \frac { \Gamma(-2s) \Gamma(-s)}{ \Gamma(1/2 +s)} ds\right| < \infty,$$
where  $C >0$ is an absolute constant.  A similar estimate holds for the repeated integral for the derivative with respect to $x$,  and the differentiation under the integral sign of the inner integral by $s$ is possible via the Stirling asymptotic formula for the gamma-function. Thus the differentiation in (2.13) is permitted as well as the change of the order of integration,  owing to the absolute and uniform convergence by $x \ge x_0 >0$,  and with the use of the reduction formula for the gamma-function we obtain 

$$\varphi(x) =   {1 \over 2\pi^3 i}  \int_{\gamma   -i\infty}^{\gamma  +i\infty} \frac { \Gamma(1 -s)}{\Gamma( 1/2 -s)}   x^{ -s} \int_0^\infty \tau \sinh(2\pi\tau)  (F_\mu f)(\tau)  \Gamma(i\tau-s) \Gamma(-i\tau -s)\  d\tau ds.\eqno(2.14)$$
Returning to integral representation (2.8) for $\varphi(x)$, we find that if the function $f^*(3/2-s)/s$ is analytic in the open strip $\gamma < {\rm Re} s < \nu+1/2$ and $(1+|s|)^{1-2{\rm Re} s} f^*(3/2-s)/ s$ is absolutely integrable over any vertical line of the closure of this strip, then via Cauchy's theorem we can move the contour to the left, integrating in (2.8) over   $(\gamma   -i\infty, \gamma  +i\infty),\  \gamma< 0 $.  Therefore  under conditions of Theorem 1 and via the above estimates we observe that both sides of (2.14) are inverse Mellin transforms of absolutely integrable functions. Hence the  uniqueness theorem can be applied for the Mellin transform in $L_1$ \cite{tit},  and we derive

$$  f^*(3/2-s) =    {1 \over \pi^2}   \frac { \Gamma(s+1/2) \Gamma (1/2+s-\mu)}{ \Gamma(1/2 -\mu-s) \Gamma(s)} $$

$$\times \int_0^\infty \tau \sinh(2\pi\tau)  (F_\mu f)(\tau) \Gamma(i\tau-s) \Gamma(-i\tau -s) d\tau,\   \hbox{max} 
\left(- {1\over 2},\ {\rm Re \mu}- {1\over 2}\right) < \gamma < 0. \eqno(2.15)$$
After the simple change of variable the latter equality becomes

$$  f^*(s) =    {1 \over \pi^2}   \frac { \Gamma(2-s) \Gamma (2- s-\mu)}{ \Gamma(s-1 -\mu) \Gamma(3/2-s)} \int_0^\infty \tau \sinh(2\pi\tau)  (F_\mu f)(\tau) \Gamma(s-3/2+ i\tau) \Gamma(s-3/2-i\tau) d\tau,$$
where   $ {\rm Re}\  s = \gamma_0= 3/2-\gamma.$ 	Hence, under condition $f \in L_{1/2-\gamma,1} (\mathbb{R}_+)$ we take  the inverse Mellin transform of both sides of the latter equality and change  the order of integration in the right-hand side of the obtained equality.  Thus we deduce,  finally,  the inversion formula for the index transform (1.4), which can be written (with the use of the reduction formula for the gamma-function and  differentiation under the integral sign ) in the form  

$$f(x)=  {x^{-1/2}  \over  \pi^2}   {d\over dx}  \  x^{3/2} \int_0^\infty \tau \sinh(2\pi\tau) S_\mu(x,\tau)  (F_\mu f)(\tau)d\tau,\ x >0,\eqno(2.16)$$
where

$$S_\mu(x,\tau)= {1 \over 2\pi i}    \int_{\gamma_0 -i\infty}^{\gamma_0+i\infty}\frac { \Gamma(2-s) \Gamma (2- s-\mu)}{ \Gamma(s-1 -\mu) \Gamma(5/2-s)}  \Gamma(s-3/2+ i\tau) \Gamma(s-3/2-i\tau)  x^{-s} ds,\ x >0.\eqno(2.17)$$
This kernel can be calculated in terms of the derivative of the product of the first and second kind associated Legendre functions. We will do it with the aid of relation (8.4.42.34) in \cite{prud}, Vol. III.  In fact, employing the  reflection formula for the gamma-function \cite{erd}, Vol. 1, we obtain 

$$  {2\sqrt \pi  e^{-i\mu\pi}\over x^2} { \Gamma(1/2-\mu+i\tau) \over \Gamma(1/2+\mu+i\tau)} P^\mu_{-1/2+i\tau} \left(\sqrt{1+x\over x}\right) \left[  Q^\mu_{-1/2+i\tau} \left(\sqrt{1+x\over x}\right) +  Q^\mu_{-1/2- i\tau} \left(\sqrt{1+x\over x}\right)\right] $$

$$=  {1 \over 2\pi i}    \int_{\gamma_0 -i\infty}^{\gamma_0+i\infty}\frac {\Gamma(2-s)\Gamma(2-\mu-s) \Gamma(s-3/2) }{ \Gamma(s-\mu-1) }  \left[ {\Gamma(s-3/2+ i\tau) \over \Gamma(5/2+i\tau-s)}  + {\Gamma(s-3/2- i\tau) \over \Gamma(5/2-i\tau-s)}\right]  x^{-s} ds$$

$$= {\cosh(\pi\tau)  \over \pi i}    \int_{\gamma_0 -i\infty}^{\gamma_0+i\infty}\frac { \Gamma(2-s) \Gamma(2-\mu-s) \Gamma(s-3/2+ i\tau) \Gamma(s-3/2-i\tau)}{ \Gamma(s-\mu-1) \Gamma(5/2-s)}   x^{-s} ds,\ x >0.$$
Hence owing to the identity \cite{erd}, Vol. 1

$$  e^{-i\mu\pi} { \Gamma(1/2-\mu+i\tau) \over \Gamma(1/2+\mu+i\tau)}  Q^\mu_{-1/2+i\tau} \left(\sqrt{1+x\over x}\right) =  e^{i\mu\pi} Q^{-\mu} _{-1/2+i\tau} \left(\sqrt{1+x\over x}\right),$$
we end up  with the following value of the kernel (2.17) 

$$S_\mu(x,\tau)=   {\sqrt \pi \  e^{i\mu\pi}\over x^2\  \cosh(\pi\tau) }  P^\mu_{-1/2+i\tau} \left(\sqrt{1+x\over x}\right) $$

$$\times\left[  Q^{-\mu}_{-1/2+i\tau} \left(\sqrt{1+x\over x}\right) +  Q^{-\mu}_{-1/2- i\tau} \left(\sqrt{1+x\over x}\right)\right],\ x > 0.$$

We summarize the results of this section as the following integral  theorem. 

{\bf Theorem 2}. {\it Let conditions of Theorem 1 hold and $\varphi \in  L_{-1/4,1}(0,1) \cap L_{1/4,1}(1,\infty),\  F_\mu  \in L_1(\mathbb{R};   |\tau| \\  \times e^{2\pi|\tau|} d\tau)$ .  Let besides,  the Mellin transform $f^*(3/2-s)$ be  such that $f*(3/2-s)/s$  be analytic in the open vertical strip  $\gamma < {\rm Re} s < \nu+1/2$ for some $\gamma <0$,  $f*(3/2-s)/s  \in L_1( ( {\rm Re} s -i\infty, {\rm Re} s +i\infty); (1+|s|)^{1-2{\rm Re} s} |ds| )$ over any vertical line of the closure  of this strip and $f \in L_{1/2-\gamma,1} (\mathbb{R}_+)$.  Then for $x >0$ the antiderivative of the function $ x^{1/2} f(x)$ can be expanded in term of the repeated integral 

$$\int_x^\infty  y^{1/2} f(y) dy = -\   { 2  e^{i\mu\pi} \over \pi \sqrt { x }}  \int_0^\infty \tau \sinh(\pi\tau) \Gamma\left({1\over 2} + i\tau -\mu\right) \Gamma\left({1\over 2} - i\tau -\mu\right) P^\mu_{-1/2+i\tau} \left(\sqrt{1+x\over x}\right) $$

$$\times\left[  Q^{-\mu}_{-1/2+i\tau} \left(\sqrt{1+x\over x}\right) +  Q^{-\mu}_{-1/2- i\tau} \left(\sqrt{1+x\over x}\right)\right] 
 \int_0^\infty \left[ P^\mu_{-1/2+i\tau} \left(\sqrt{1+y\over y}\right)\right]^2 f(y) dy\  d\tau,\eqno(2.18)$$
generating the pair $(1.4), (2.16)$ of the direct and inverse index transforms, respectively. }

{\bf Remark 1}. {\it The classical Lebedev expansion $(1.1)$ for  adjoint kernels  can be obtained, letting $\mu=0$ in (2.18).  Then,    substituting    $\sqrt{(1+x)/ x},\  \sqrt{(1+y)/y}$ by new variables, we derive

$$\int_1^x h(y)dy = 2(x^2-1)^{1/2} \int_0^\infty \tau \tanh(\pi\tau) P_{-1/2+i\tau} (x) \left[  Q_{-1/2+i\tau} (x) +  Q_{-1/2- i\tau} (x)\right]  $$

$$\times \int_1^\infty (y^2-1)^{1/2}  \left[ P_{-1/2+i\tau} (y)\right]^2  h(y) dy d\tau,$$
where $h(y)= 2y  f\left( (y^2-1)^{-1}\right) (y^2-1)^{-5/2}$. }

{\bf Remark 2}. {\it By the same technique one can establish the generalized Lebedev   expansion $(1.1)$ in the form

$$\int_x^\infty y^{1/2}  f(y)dy =    -\   { 2  e^{i\mu\pi} \over \pi \sqrt { x }}   \int_0^\infty \tau \sinh (\pi\tau)\Gamma\left({1\over 2} + i\tau -\mu\right) \Gamma\left({1\over 2} - i\tau -\mu\right) \left[ P^\mu_{-1/2+i\tau} \left(\sqrt{1+x\over x}\right)\right]^2 $$

$$\times \int_0^\infty  P^\mu_{-1/2+i\tau} \left( \sqrt{1+y\over y} \right) \left[  Q^{-\mu}_{-1/2+i\tau} \left(\sqrt{1+y\over y}\right) +  Q^{-\mu}_{-1/2- i\tau} \left(\sqrt{1+y\over y}\right)\right] f(y) dy d\tau,\ x > 0,$$
leading to the reciprocal formulas $(1.2), (1.3)$  when $\mu=0$. We leave details to the reader. }

\section{Index transform (1.5)} 

The boundedness of the adjoint operator (1.5) and inversion formula for this transformation will be established below.    We begin with 

{\bf Theorem 3.}  {\it Let ${\rm Re}\mu < 1/2$. The index transform  $(1.5)$  is well-defined as a  bounded operator $G_\mu: L_{1} \left(\mathbb{R}\right) \to L_{\nu,\infty}  (\mathbb{R}_+),\  \nu \in (-1/2,\ - {\rm Re}\mu)$ and the following norm inequality takes place
$$||G_\mu g||_{\nu,\infty}  \equiv \hbox{ess sup}_{x >0} | x^\nu (G_\mu g)(x)| \le C_{\mu,\nu}  ||g ||_{L_1(\mathbb{R})},\eqno(3.1)$$
where $C_{\mu,\nu}$ is defined by $(2.6)$.   Moreover,   if $(G_\mu g) (x) \in L_{\nu,1}(\mathbb{R}_+)$, then for all $y >0$ }
$${1\over 2\pi i}   \int_{\nu -i\infty}^{\nu  +i\infty} \frac{\Gamma(1+s) \Gamma (1+s-\mu) \Gamma(-s)} { \Gamma(1/2+s)\Gamma(-\mu-s) } (G^*_\mu g) (s) y^{ -s} ds\
=  \sqrt {\pi y}  \   \int_{-\infty}^\infty  e^{y/2} \   K_{i\tau} \left({y\over 2} \right)  {g(\tau) \over \cosh(\pi\tau) } d\tau. \eqno(3.2)$$

\begin{proof}   The norm inequality (3.1) is a direct consequence of the estimates (2.9). Indeed, since the kernel (1.6) has a bound

$$|\Phi_\tau (x)| \le C_{\mu,\nu} x^{-\nu},\ x >0,$$
where $C_{\mu,\nu}$ is defined by (2.6), we have 

$$|(G_\mu g)(x)| \le C_{\mu,\nu} x^{-\nu} \int_\mathbb{R} |g(\tau)| d\tau =  C_{\mu,\nu} x^{-\nu} ||g ||_{L_1(\mathbb{R})},$$
and (3.1) follows.  Now, taking the Mellin transform (2.1) of both sides of (1.5), we change the order of integration by Fubini's theorem and take into account (1.7) to obtain the equality

$$  \frac{\Gamma(1+s) \Gamma (1+s-\mu) } { \Gamma(1/2+s)\Gamma(-\mu-s) } (G^*_\mu g) (s)=   \int_{-\infty}^\infty \Gamma\left(s+ {1\over 2} + i\tau\right) \Gamma\left(s+ {1\over 2} - i\tau \right)  g(\tau) d\tau,\eqno(3.3)$$
where $-1/2 < {\rm Re} s < - {\rm Re} \mu$.    Hence the inverse Mellin transform (2.3) and  relation (8.4.23.5) in \cite{prud}, Vol. III will lead us to (3.2), completing the proof of Theorem 3. 
\end{proof} 

The inversion formula for the index transform (1.5) is given by

{\bf Theorem 4}.  {\it Let ${\rm Re}\mu < 1/2$,  $g(z/i)$ be an even analytic function in the strip $D= \left\{ z \in \mathbb{C}: \ |{\rm Re} z | < \alpha < 1/2\right\}, \\  such\  that  g(0)=g^\prime (0)=0$ and $g(z/i)$ be  absolutely  integrable over any vertical line in  $\overline{D}$.   If $(G_\mu g) (y) \in  L_{1- \gamma, 1} \left(\mathbb{R}_+\right),\   3/2  <  \gamma  <  \hbox{min} (2,\   2 - {\rm Re} \mu)$, then for all  $x \in \mathbb{R}$ the  inversion formula holds for the index transform (1.5)} 
$$ g(x)  = {e^{i\mu\pi} \over \pi\sqrt \pi}  \   x \sinh(\pi x) \int_0^\infty \left( y^{-1/2} \ {d\over dy}\  y^{-1/2}\right) \left[  P^\mu_{-1/2+ix} \left(\sqrt{1+y\over y}\right) \right.$$

$$\times\left. \left[  Q^{-\mu}_{-1/2+ix} \left(\sqrt{1+y\over y}\right) +  Q^{-\mu}_{-1/2- ix} \left(\sqrt{1+y\over y}\right)\right] \right] (G_\mu g)(y) dy.\eqno(3.4)$$

\begin{proof}    Indeed,  recalling  (3.2), we multiply its both sides by $ e^{-y/2} K_{ix} \left({y/2} \right) y^{\varepsilon - 3/2}$ for some positive $\varepsilon \in (0,1)$ and integrate with respect to $y$ over $(0, \infty)$.  Hence changing the order of integration in the left-hand side of the obtained equality due to the absolute convergence of the  iterated integral, we  appeal  to relation (8.4.23.3) in \cite{prud}, Vol. III to find  
$$ {1\over 2\pi i}   \int_{\nu -i\infty}^{\nu  +i\infty} \frac { \Gamma(1+s) \Gamma (1+s-\mu) \Gamma(\varepsilon- s-1/2 + ix)\Gamma(\varepsilon- s-1/2 -ix) \Gamma(-s)}{\Gamma(1/2+s)  \Gamma( \varepsilon -s) \Gamma(-\mu-s) } (G_\mu g)^*(s) ds$$
$$ =  \ \int_0^\infty  K_{ix} \left({y\over 2} \right) y^{\varepsilon -1} \int_{-\infty}^\infty  K_{i\tau} \left({y\over 2} \right)  {g(\tau) \over \cosh(\pi\tau) } d\tau dy. \eqno(3.5)$$
In the meantime,  the right-hand side of (3.5) can be treated, using   the evenness of $g$ and  the representation  of the modified Bessel  function $K_z(y)$ in terms of the modified Bessel function of the first kind $I_z(y)$ \cite{erd}, Vol. II. Hence with a simple substitution we find 
$$  \int_0^\infty  K_{ix} \left({y\over 2} \right) y^{\varepsilon -1} \int_{-\infty}^\infty  K_{i\tau} \left({y\over 2} \right)  {g(\tau) \over \cosh(\pi\tau) } d\tau dy$$

$$= 2 \pi i \int_0^\infty  K_{ix} \left({y\over 2} \right) y^{\varepsilon -1} \int_{-i\infty}^{i\infty}   I_{ z} \left({y\over 2} \right)  {g(z/i) \over \sin (2\pi z) } dz\  dy. \eqno(3.6)$$
On the other hand, according to our assumption $g(z/i)$ is analytic in the vertical  strip $0\le  {\rm Re}  z < \alpha< 1/2$,  $g(0)=g^\prime (0)=0$ and integrable  in the closure of the strip.  Hence,  appealing to the inequality for the modified Bessel   function of the first  kind  (see \cite{yal}, p. 93)
 $$|I_z(y)| \le I_{  {\rm Re} z} (y) \  e^{\pi |{\rm Im} z|/2},\   0< {\rm Re} z < \alpha,$$
one can move the contour to the right in the latter integral in (3.6). Then 

$$2 \pi i  \int_0^\infty  K_{ix} \left({y\over 2} \right) y^{\varepsilon -1} \int_{-i\infty}^{i\infty}   I_{ z} \left({y\over 2} \right)  {g(z/i) \over \sin (2\pi z) } dz\  dy$$

$$= 2 \pi i    \int_0^\infty  K_{ix} \left({y\over 2} \right) y^{\varepsilon -1} \int_{\alpha -i\infty}^{\alpha + i\infty}   I_{ z} \left({y\over 2} \right)  {g(z/i) \over \sin (2\pi z) } dz\  dy.$$
Now ${\rm Re} z >0$,  and  it is possible to pass to the limit under the integral sign when $\varepsilon \to 0$ and to change the order of integration due to the absolute and uniform convergence.  Therefore the value of the integral (see relation (2.16.28.3) in \cite{prud}, Vol. II)
$$\int_0^\infty K_{ix}(y) I_z(y) {dy\over y} = {1\over x^2 + z^2} $$ 
leads us to the equalities 

$$\lim_{\varepsilon \to 0}  2 \pi i  \int_0^\infty  K_{ix} \left({y\over 2} \right) y^{\varepsilon -1} \int_{-i\infty}^{i\infty}   I_{ z} \left({y\over 2} \right)  {g(z/i) \over \sin (2\pi z) } dz\  dy$$

$$=    2 \pi i  \int_{\alpha -i\infty}^{\alpha + i\infty}   {g(z/i) \over (x^2+ z^2) \sin (2\pi z) } dz =  \pi   i 
 \left( \int_{-\alpha +i\infty}^{- \alpha- i\infty}   +   \int_{\alpha -i\infty}^{ \alpha+  i\infty}   \right)  {  g(z/i) \  dz \over (z-ix) \  z \sin(2\pi z)}. \eqno(3.7)$$
Hence conditions of the theorem allow to apply the Cauchy formula in the right-hand side of the latter equality in (3.7).  Thus 
$$\lim_{\varepsilon \to 0}  2 \pi i   \ \int_0^\infty  K_{ix} \left({y\over 2} \right) y^{\varepsilon -1} \int_{-i\infty}^{i\infty}   I_{ z} \left({y\over 2} \right)  {g(z/i) \over \sin (2\pi z) } dz\  dy =  { 2\pi^{2} \  g(x) \over  x\sinh (2\pi x)} ,\quad x \in \mathbb{R} \backslash \{0\}.\eqno(3.8)$$
On the other hand,  appealing to  the Parseval identity  (2.2), the left-hand side of (3.5) can be rewritten in the form
$$ {1\over 2\pi i}   \int_{\nu -i\infty}^{\nu  +i\infty} \frac { \Gamma(1+s) \Gamma (1+s-\mu) \Gamma(\varepsilon- s-1/2 + ix)\Gamma(\varepsilon- s-1/2 -ix) \Gamma(-s)}{\Gamma(1/2+s)  \Gamma( \varepsilon -s) \Gamma(-\mu-s) } (G_\mu g)^*(s) ds$$

$$ = \int_0^\infty (G_\mu g)(y) \Psi_{\varepsilon}  (x,y)  dy,\eqno(3.9)$$
where
$$\Psi_{\varepsilon}  (x,y) = {1\over 2\pi i}   \int_{1-\nu -i\infty}^{1-\nu  +i\infty} \frac {\Gamma(\varepsilon+ s-3/2 + ix)\Gamma(\varepsilon+ s-3/2 -ix) \Gamma(s-1)  \Gamma(2-s) \Gamma (2-s-\mu) }{\Gamma(s+ \varepsilon -1) \Gamma(3/2-s)   \Gamma(s-\mu-1) } y^{-s} ds.\eqno(3.10)$$ 
However,  due to the properties of the integrand in (3.10) as an analytic and absolutely integrable function in some strip,  one can move the contour to the right within the strip $\hbox{max} (1, \  3/2- \varepsilon) < {\rm Re} s <  \hbox{min} (2,\   2 - {\rm Re} \mu)$.  This circumstance allows to calculate the limit when $\varepsilon \to 0+$ in (3.10), passing under the integral sign.  Hence, comparing with (2.17), we find

$$\lim_{\varepsilon \to 0+} \Psi_{\varepsilon}  (x,y) =  {\sqrt \pi \  e^{i\mu\pi}\over \  \cosh(\pi x) }   \left( y^{-1/2} \ {d\over dy}\  y^{-1/2}\right) \left[  P^\mu_{-1/2+ix} \left(\sqrt{1+y\over y}\right) \right.$$

$$\times\left. \left[  Q^{-\mu}_{-1/2+ix} \left(\sqrt{1+y\over y}\right) +  Q^{-\mu}_{-1/2- ix} \left(\sqrt{1+y\over y}\right)\right] \right].$$
Further,  the passage to the limit under integral sign when $\varepsilon \to 0+$   in the right-hand side of (3.9) is possible, appealing to the dominated  convergence theorem. In fact,  due to the estimate with the use of the definition of Euler's beta-function and the duplication formula for the gamma-function, we have 

$$\left| \Psi_{\varepsilon}  (x,y) \right| \le {y^{- {\rm Re} s} \  2^{ 2 {\rm Re} s -3}  \over \pi \sqrt \pi }   B\left( {\rm Re} s - 3/2,  \  {\rm Re} s - 3/2\right) \Gamma( \varepsilon+  {\rm Re} s-3/2)$$

$$\times   \int_{{\rm Re} s -i\infty}^{{\rm Re} s  +i\infty} \left| \frac { \Gamma(s-1)  \Gamma(2-s) \Gamma (2-s-\mu) }{ \Gamma(3/2-s)   \Gamma(s-\mu-1) }  ds\right|\le C \  y^{- {\rm Re} s},$$
where $C >0$ is an absolute constant,   $3/2  < {\rm Re} s =\gamma  <  \hbox{min} (2,\   2 - {\rm Re} \mu)$.  Together with  the condition $(G_\mu g)(y) \in L_{1- \gamma, 1} \left(\mathbb{R}_+\right)$  it gives the desired property.   Finally, combining with (3.8),  (3.5), we establish the inversion formula (3.4). Theorem 4 is proved. 

 \end{proof}

\section{Boundary   value problem}

In this section  the index transform (1.5) is employed to investigate  the  solvability  of the boundary  value  problem  for the  following third   order partial differential  equation

$$ x \left( r(1+r) - y^2\right)  {\partial^3 u \over \partial x^3} +   y \left( r(1+r) - x^2\right)  {\partial^3 u \over \partial y^3}$$

$$ +   x \left( r(1+r) -  x^2+ 2y^2\right)  \ {\partial^3 u \over \partial x\partial y^2} +  y \left( r(1+r)  - y^2 + 2x^2\right)   {\partial^3 u \over \partial y \partial x^2}$$

$$+   \left( {3\over 2 r} + 2 \right)  \left( \left( y^2+ {x^2\over r}\right) \  {\partial^2 u \over \partial x^2} + \left( x^2+ {y^2\over r}\right)\   {\partial^2 u \over \partial y^2} \right)$$

$$  + \  {(3+r)\  xy  \over r^2}  \    {\partial^2 u \over \partial x\partial y}  - \left(  \mu^2 -1 +  {1\over 4 r} \right)  \left( x {\partial u \over \partial x}  + y {\partial u  \over \partial y}\right)   - {u\over 8r}  =0,\eqno(4.1)$$
where $(x,y) \in \mathbb{R}^2 \backslash \{ 0\},\  r= \sqrt {x^2+y^2}$.    Writing  (4.1) in polar coordinates $(r,\theta)$, we end up with the equation

$$r^2(1+r)  {\partial^3  u \over \partial r^3} +  \  {\partial^3  u \over \partial r \partial \theta^2} +  3 \left(r+ {1\over 2}\right) {\partial^2 u  \over \partial r^2}  -  {1\over 2 r}\  {\partial^2  u \over  \partial \theta^2} +   \left( r(1-\mu^2) +  {1\over 4} \right)  {\partial u \over \partial r}  -  {u\over 8r }  = 0.\eqno(4.2)$$

{\bf Lemma 3.} {\it  Let ${\rm Re}\  \mu < 0,\  g(\tau)  \in L_1\left(\mathbb{R}; \  e^{ \beta  |\tau|} d\tau\right),\  \beta \in (0, \pi)$. Then  the function
$$u(r,\theta)=    \sqrt\pi \int_{-\infty}^\infty \Gamma\left({1\over 2} + i\tau -\mu\right) \Gamma\left({1\over 2} - i\tau -\mu\right)  \left[ P^\mu_{-1/2+i\tau} \left(\sqrt{1+r\over r}\right)\right]^2 e^{\theta \tau} \ g(\tau) d\tau,\eqno(4.3)$$
 satisfies   the partial  differential  equation $(4.2)$ on the wedge  $(r,\theta): r   >0, \  0\le \theta <  \beta$, vanishing at infinity.}

\begin{proof} The proof  is straightforward by  substitution (4.3) into (4.2) and the use of  (1.12).  The necessary  differentiation  with respect to $r$ and $\theta$ under the integral sign is allowed via the absolute and uniform convergence, which can be verified, appealing to  the integrability condition $g \in L_1\left(\mathbb{R}; e^{ \beta |\tau|} d\tau\right),\  \beta \in (0, \pi)$ and estimates of the derivatives of the kernel (1.6) with respect to $r$. In fact,  it is  based on bounds of the corresponding integral of type (1.7)  with the use of equalities (1.10),  (1.11).   Finally,  the condition $ u(r,\theta) \to 0,\ r \to \infty$  is due to the asymptotic behavior of the associated Legendre function of the first kind near  unity in the case  ${\rm Re}\  \mu < 0$ (see \cite{erd}, Vol. I.)  
\end{proof}

Finally  we will formulate the boundary  value problem for equation (4.2) and give its solution.

{\bf Theorem 5.} {\it Let  $g(x)$ be given by formula $(3.4)$ and its transform $(G_\mu g) (y)\equiv G_\mu (y)$ satisfies conditions of Theorem 4.  Then  $u (r,\theta),\   r >0,  \  0\le \theta < \beta$ by formula $(4.3)$  will be a solution  of the boundary  value problem on the wedge for the partial differential  equation $(4.2)$ subject to the boundary  condition}
$$u(r,0) = G_\mu(r).$$

\bigskip
\centerline{{\bf Acknowledgments}}
\bigskip

The work was partially supported by CMUP (UID/MAT/00144/2013), which is funded by FCT(Portugal) with national (MEC) and European structural funds through the programs FEDER, under the partnership agreement PT2020.

\bibliographystyle{amsplain}

\begin{thebibliography}{10}

\bibitem{square}  Lebedev  NN.   Expansion of an arbitrary function in an integral with respect to the squares of Legendre functions with complex index,  {\it Differencial'nye uravneniya},   {\bf 3}  (1967),  422-435  (in Russian).

\bibitem{yak}   Yakubovich S.  Index transforms.   Singapore:  World Scientific Publishing Company; 1996.

\bibitem{vir}   Virchenko N, Fedotova I.  Generalized associated Legendre functions and their applications (with a Foreword by Semyon Yakubovich).   Singapore:  World Scientific Publishing Company;  2001.

\bibitem{erd}    Erd\'elyi A,  Magnus W,   Oberhettinger  F,   Tricomi FG.  Higher transcendental functions. Vols. I,  II. New  York: McGraw-Hill;  1953.

\bibitem{prud}  Prudnikov AP,  Brychkov  YuA,  Marichev OI. Integrals and series:  Vol. I: Elementary functions. New York:  Gordon and Breach;   1986;   Vol. II:  Special functions. New York: Gordon and Breach;  1986;   Vol. III:  More special functions. New York:   Gordon and Breach; 1990.

\bibitem{yal}   Yakubovich S,  Luchko Yu.  The hypergeometric approach to integral transforms and convolutions, Mathematics and its applications.  Vol. 287.  Dordrecht:  Kluwer Academic Publishers Group; 1994.


\bibitem {tit}  Titchmarsh EC.   An introduction to the theory of Fourier integrals.   New York:  Chelsea; 1986.

\bibitem{square}  Lebedev N.N.   On an integral representation of an arbitrary function in terms of squares of Macdonald functions with imaginary index,  {\it Sibirsk. Mat. Zh.},   {\bf 3}  (1962),  213-222 (in Russian).










\end{thebibliography}

\end{document}